\documentclass{article}

\usepackage{maa-monthly-modified}
\theoremstyle{theorem}
\newtheorem{theorem}{Theorem}
\newtheorem{corollary}{Corollary}
\newtheorem{proposition}{Proposition}
\newtheorem{lemma}{Lemma}

\theoremstyle{definition}

\begin{document}

\title{Navigating Around Convex Sets}
\markright{Navigating Convex Sets}
\author{J. J. P. Veerman}

\maketitle

\noindent
\section*{Abstract}
\begin{small}
We review some basic results of convex analysis and geometry in $\mathbb{R}^n$ in the context of
formulating a differential equation to track the distance between an observer flying outside
a convex set $K$ and $K$ itself.
\end{small}

%\vskip 0.4in\noindent
\section{Introduction.}
\label{chap:intro}
\setcounter{figure}{0} \setcounter{equation}{0}

\noindent
Suppose you fly by a convex body in $\mathbb{R}^n$. A fundamental problem in convex analysis is to
determine your distance $r(t)$ to that body as you are flying.
To the best of our knowledge, all standard algorithms to do this determine the point in the convex set
nearest to you by optimization, at every point in time when needed or possible. In this exposition,
we propose to solve the optimization problem only once, at $t=0$, say, and then track $r(t)$
by solving a differential equation for it. The study of this problem allows us to
review some of the basic tenets of convex analysis.

Let $K$ be a closed (solid) convex body in $\mathbb{R}^n$ and denote by $\mathbb{R}^n\backslash K$
by $\Omega$. The fly-by curve in $\Omega$ is given by the trajectory $c(t)$, which we will always assume
to be smooth. To find $r(t)$, we need to know the point in $K$ nearest to $c(t)$.
This point is called the \emph{projection of $c(t)$ onto $K$} and will be denoted by $\Pi(c(t))$ or $\Pi(t)$
for short. If we hope to write down a differential equation for $r(t)$, we will \emph{at least} have
to know the one-sided derivative of $\Pi(c(t))$ with respect to time:
\begin{equation}
\Pi_+'(t_0)\equiv \,\lim_{t\searrow 0}\, \dfrac{\Pi(t)-\Pi(0)}{t} .
\label{eq:one-sidedetc}
\end{equation}
In Section \ref{chap:C} we give an elementary proof that if in $\mathbb{R}^2$, the boundary
$\partial K$ of $K$ is twice differentiable, then this one-sided derivative exists.
In Section \ref{chap:counter-example}, we will see that this derivative exists for
any piecewise linear polygon. Then we will outline the construction from \cite{AMV} of
a convex body
in $\mathbb{R}^2$ whose boundary is $C^{1,1}$ ($C^1$ with a Lipschitz derivative) but for which
$\Pi_+'$ does not exist (Theorem \ref{thm:no-directional-deriv}).

In Sections \ref{chap:R2} and \ref{chap:R3}, the differential equations are constructed
that allow us to continually monitor the distance $r(t)$ to a convex body as we navigate around it.
The general case is not much different from $\mathbb{R}^3$, so for ease of exposition we stick
to $\mathbb{R}^3$. The construction is straightforward and employs the Weingarten equations
(or Ricci curvature) from differential geometry. To the best of our knowledge, the final form of
these equations (Theorem \ref{thm:3dim}) is new. We invoke the existence and uniqueness theorems for
the solutions of differential equations to show that if the boundary of $K$ is $C^{2,1}$, then
the system of differential equations has a unique solution. In Section \ref{chap:examples}, we
give two examples of these equations and their solutions.

To the careful reader it will be clear that all we really need to formulate these differential
equations are one-sided time derivatives (as in equation (\ref{eq:one-sidedetc})).
Thus if you fly near a convex body, you can change your course to another differentiable one,
all the while continuing to monitor the distance. Thus the equations of Sections \ref{chap:R2}
and \ref{chap:R3} can be used to \emph{avoid} collision with a convex set. Hence our title.
Even if $K$ is not convex, one could use these equations to avoid the convex hull of $K$.

Finally, in Section \ref{chap:dist-diffb}, we discuss the classical, but still remarkable, fact
that the function $r:\Omega \rightarrow \mathbb{R}$ giving the distance of a point to the \emph{convex}
set is \emph{always} differentiable (Theorem \ref{thm:differentiable})!
There is no regularity requirement at all on the convex set. Our exposition of this fact is inspired
by \cite{foote}.

The subject of distance to convex sets is much too rich to do justice to in this short
introduction. But we do wish to point out a few directions in which considerable further
research has been done. Our own interest (beside \cite{AMV}) ultimately derives from a
closely related problem, namely: Given a Riemannian manifold and two disjoint compact
subsets $A$ and $B$, what are the topological \cite{BV,VB} and geometric \cite{HPV} characteristics
of the set of points whose distance to $A$ equals their distance to $B$?
Beside their intrinsic interest, these sets have many applications among others in the study of
Brillouin zones in quantum mechanics \cite{peix}.

The distance function to $\partial K$ is \emph{not} differentiable
\emph{at} $\partial K$. Indeed, at a point of $\partial K$, it behaves like the function $x\rightarrow |x|$
near 0. For this reason one is also interested \cite{ambrosio} in the smoothness of the \emph{signed}
distance, which is negative on one side of the surface and positive on the other. This obviously
only works for embedded codimension one surfaces. Thus for more generality, one also studies
the properties of $r^2$. This is not formally a distance, but still closely related, and
is differentiable at ${\partial K}$ (see \cite{ambrosio}).

There are many generalizations of the regularity of the distance function. In a smooth Riemannian manifold,
or Alexandrov space, the notion of convexity may not be well-defined. So, instead of the derivative of the distance
$r(t)$ to a closed subset $K$, one looks instead at its one-sided derivative along a geodesic $c(t)$
with initial point $x_0\in \Omega$:
\begin{equation}
r_+'(0)\equiv \,\lim_{t\searrow 0}\, \dfrac{r(c(t))-r(c(0))}{t} .
\label{eq:one-sided2}
\end{equation}
The generalization of Theorem \ref{thm:differentiable} holds and states that this one-sided derivative
always exists (see \cite{Plaut} and \cite[Exercise 4.5.11]{BBI}). Further differentiability beyond that
depends on the smoothness of the subset $K$ (for example, \cite{krantz,nirenberg}). Other
generalizations consider the distance to convex sets in Hilbert spaces (see \cite{fitzpatrick,shapiro2}).

%\vskip 0.4in\noindent
\section{Twice Differentiable Sets in $\mathbb{C}$.}
\label{chap:C}

\begin{figure}[pbth]
\begin{center}
\includegraphics[width=5.0cm]{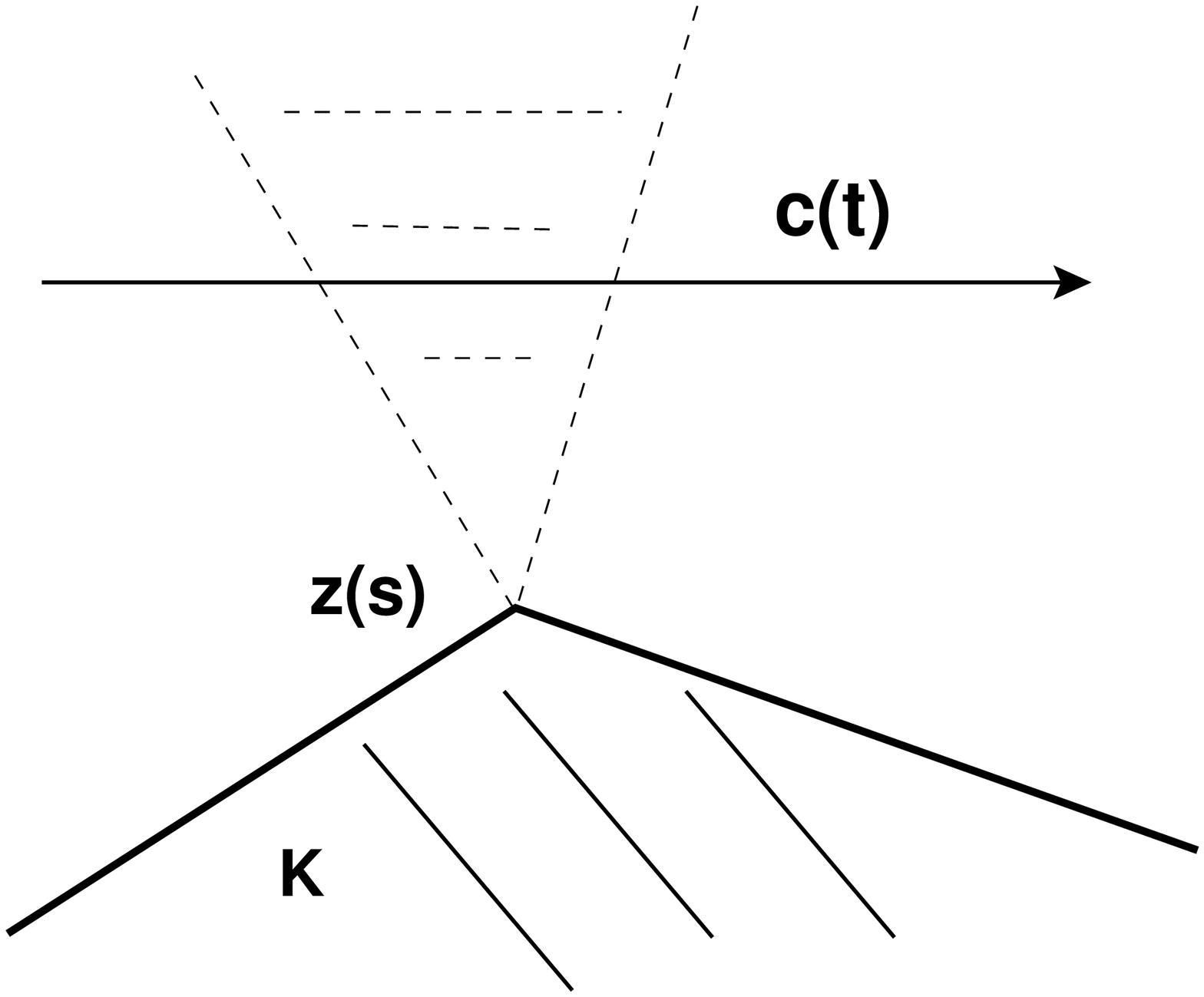} \quad \quad \quad \quad
\includegraphics[width=3.3cm]{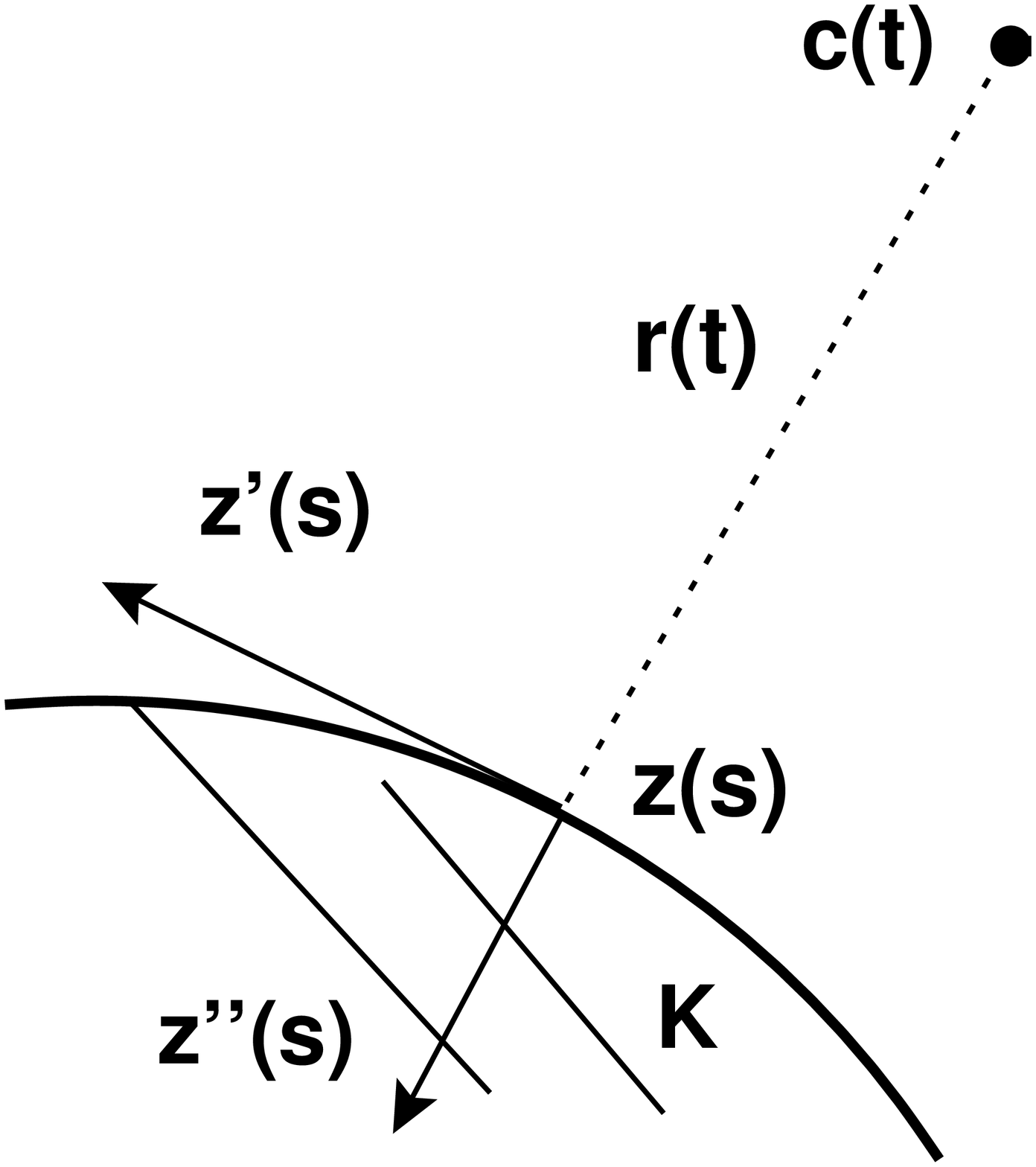}
\caption{In the left figure, $K$ has a nondifferentiable boundary, but the one-sided derivative of the
projection still exists. In the right figure, $K$ has a twice differentiable boundary.
The tangent vector $z'(s)$
and the curvature vector $z''(s)$ are orthogonal, because $z$ is a unit speed parametrization.}
\label{fig:convexset}
\end{center}
\end{figure}

\noindent
It is convenient to use complex coordinates in $\mathbb{R}^2$. So, we identify $\mathbb{R}^2$ with $\mathbb{C}$. We will assume that the boundary $\partial K$ is a twice differentiable curve $z(s)$.
We also assume that $z:\mathbb{R}\rightarrow \mathbb{C}$ is a unit speed parametrization, that is,
$|z'(s)|=1$. We orient the boundary $z(s)$ of $K$ \emph{counter-clockwise} (see Figure \ref{fig:convexset}).

The trajectory of a point \emph{outside the body $K$} depends on time and is given by
\begin{equation}
c=z-irz' \quad \textrm{or} \quad  c(t)=z(s(t))-ir(t)z'(s(t)) ,
\label{eq:coords}
\end{equation}
so that in $\Omega$ (outside the convex body), $r$ is positive.
Note that we indicate differentiation with respect to time $t$ with a dot (for example, $\dot c(t)$),
whereas differentiation with respect to the parameter of the convex body is indicated
by an accent (for example, $z'(s)$).
\begin{equation}
\dot c = (z'-irz'')\dot s - i \,\dot r z' .
\label{eq:cdot}
\end{equation}
Now denote by $\dot c_{\|}$ and $\dot c_{\bot}$ the components of $\dot c$ parallel and
orthogonal, respectively, to $z'$ (derivative with respect to $s$). Since $|z'(s)|=1$, we have that
$z''(s)$ is orthogonal to $z'(s)$, and thus $iz''(s)$ is parallel to $z'$. Thus, equation
(\ref{eq:cdot}) quite naturally splits into two components as described in the following lemma:
\begin{lemma} Suppose $z$ is twice differentiable
and $c$ is a differentiable trajectory outside $K$. We have
\begin{equation}
\left\{\begin{array}{ccccc}
\dot s &=& \dfrac{\dot c_{\|}}{z'-irz''}&=& \textrm{Re}\left(\dfrac{\dot c}{z'-irz''}\right)\\[0.35cm]
\dot r &=& \dfrac{i\dot c_{\bot}}{z'} &=& \textrm{Re}\left(\dfrac{i\dot c}{z'}\right)\\[0.2cm]
\end{array} \right. \quad .
\label{eq:ode}
\end{equation}
\label{lem:complex}
\end{lemma}

\vskip -0.2in

%\vskip -0.2in
Equation (\ref{eq:coords}) implies that in the orientation of the parametrization sketched in
Figure \ref{fig:convexset}, $r(t)$ is positive. Thus $irz''$ is parallel to $z'$,
but actually points in the opposite direction. Because $z'(s) \dot s$
is the derivative of the projection, equation (\ref{eq:ode}) proves the following lemma.
This result generalizes to higher dimension (see Corollary \ref{cor:deriv-proj}).

\begin{lemma} If $\partial K$ is twice differentiable, then $|\Pi_+'(0)|=\dfrac{|\dot c_{\|}|}{1+r\,|z''|}$
($r>0$ in $\Omega$).
\label{lem:curvature}
\end{lemma}

%\vskip 0.4in\noindent
\section{A $C^{1,1}$ Counterexample.}
\label{chap:counter-example}

\noindent
If we drop the requirement that $z(s)$ is twice differentiable, a set can be constructed where
$\Pi_+'(0)$ does not exist. We now give that construction, and outline the reason it works.
Our treatment is loosely based on \cite{AMV}, which contains the full details.

One might be tempted to say that, of course, $z$ must be differentiable for $\Pi_+'$ to exist.
But that would be wrong! Consider the set up on the left side of Figure \ref{fig:convexset},
where $c(t)$ is a straight line ($c'(t)$ is constant).
Consider the (shaded) cone in $\Omega$ formed by the perpendiculars to the tangents of the
surface at the nondifferentiable point. Before reaching the cone, $\Pi_+'(t)$ is a constant
determined by the angle between $c'(t)$ and the left flank of $\partial K$. As soon as
$c$ hits the left boundary of the cone, $\Pi_+'(t)=0$ because we are looking at the one-sided
derivative (see equation (\ref{eq:one-sidedetc})). Again, when $c$ hits the right boundary of the
cone, the one-sided derivative is a nonzero constant. Clearly, the one-sided derivative exists
for every piecewise linear polygon. In the light of this, it becomes nontrivial
to find a counterexample to the existence of the one-sided derivative of the projection.
A beautiful and surprisingly simple example of a nonempty closed \emph{continuous} convex set
for which the directional derivative of the metric projection mapping fails to exist was constructed
by Shapiro in \cite{shapiro} and fine-tuned in \cite{AMV} to become $C^{1,1}$.

\begin{figure}[pbth]
\begin{center}
\includegraphics[width=5.5cm]{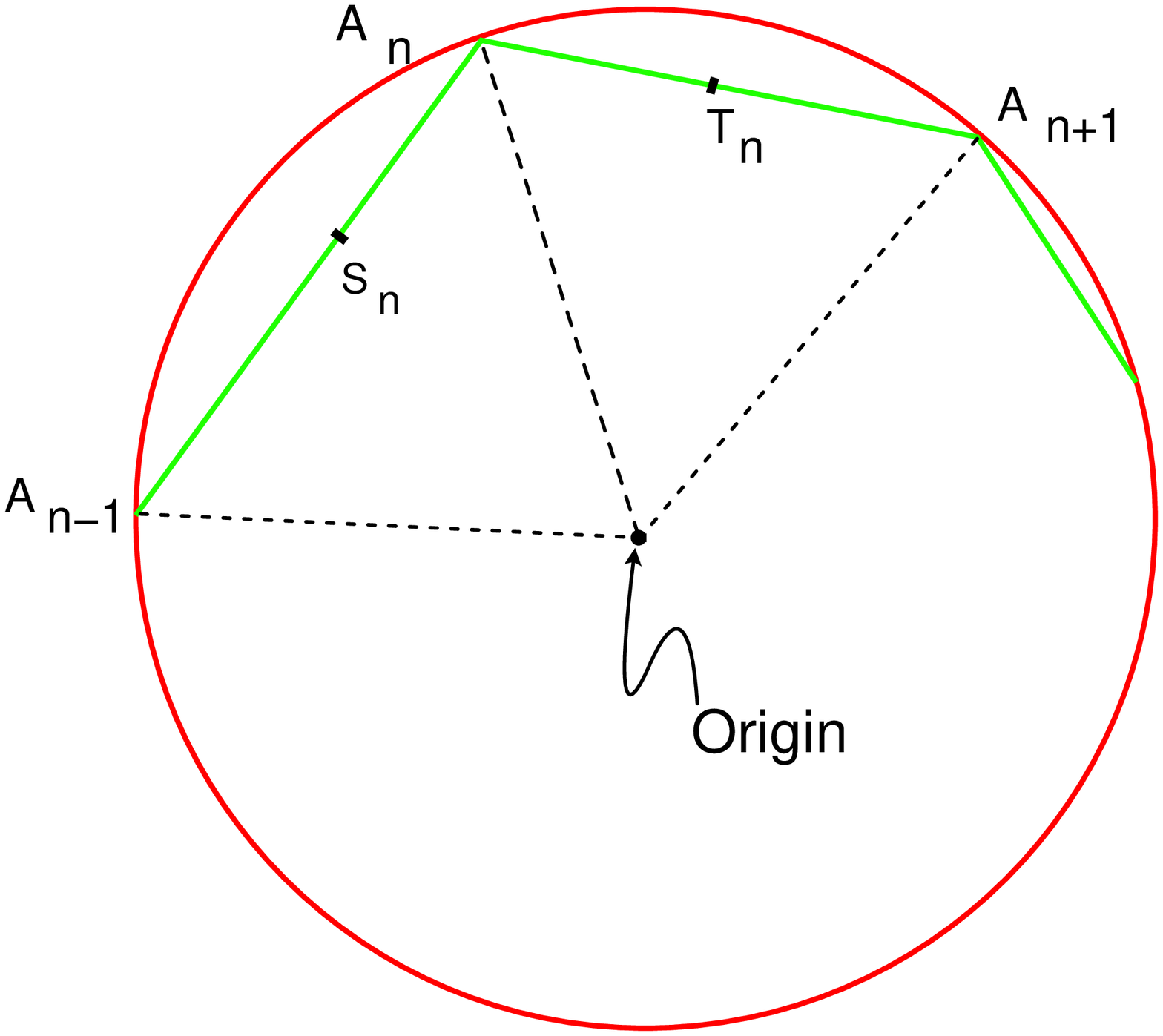} \quad \quad \quad
\includegraphics[width=5.5cm]{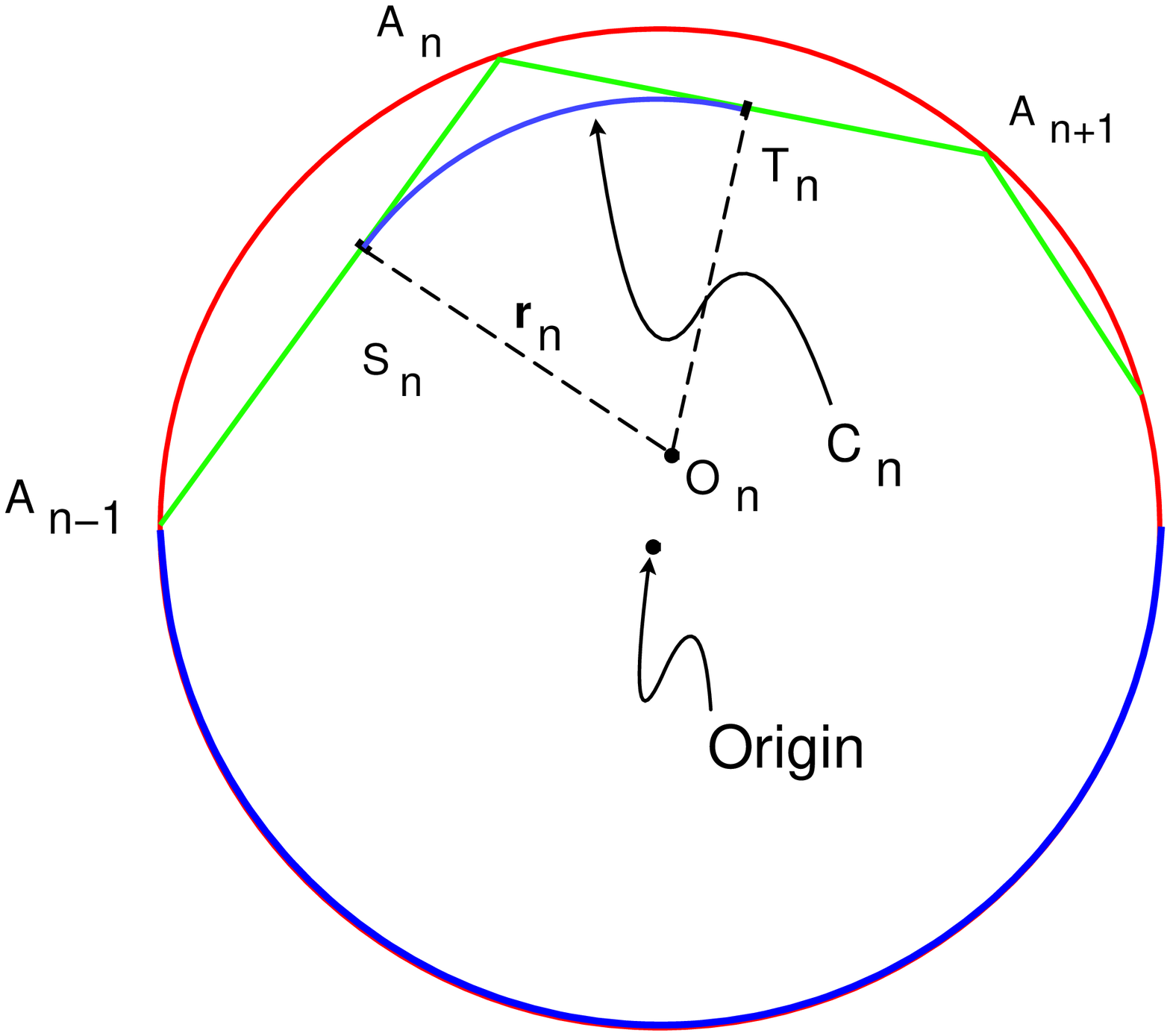}
\caption{On the left, the $C^0$ counterexample given by Shapiro, and on the right its modification
into a $C^{1,1}$ set. For clarity, the drawing is not on a realistic scale.}
\label{fig:convexset1}
\end{center}
\end{figure}

For some $\lambda\in(0,1)$, define a  sequence of real numbers
$\{\alpha_n\}_{n\in\mathbb{N}}\subset (0,\pi]$ by
\begin{equation}
\alpha_n = C\lambda^n  \quad \textrm{and set } \quad A_n=e^{i\alpha_n} .
\label{eq:alpha}
\end{equation}
($C$ is any positive real so that $\alpha_n \in (0,\pi]$ for all $n$.)
This set in \cite{shapiro} is the convex hull of the collection of points $0$, $1$, and
$\{A_n\}_{n\in\mathbb{N}}$. See the left part of Figure \ref{fig:convexset1}.

The right part of Figure \ref{fig:convexset1} shows the figure modified by the techniques of
\cite{AMV}, which is now $C^{1,1}$.
The construction is as follows. Let $T_n$ be the \emph{midpoint} of the line
segment $A_nA_{n+1}$ and let $S_n$ be the point in the line segment $A_{n-1}A_n$ such that
\begin{displaymath}
\textrm{length}(T_nA_n)= \textrm{length}(A_nS_n) .
\end{displaymath}
Replace the two line segments $T_nA_n$ and $A_nS_n$ by a circular arc $C_n$ tangent to both
segments. Define the convex set $K$ as the convex hull of all arcs $C_n$
and the arc of the unit circle in the lower half-plane connecting $1$ and $A_1$.
Clearly, the boundary of $K$ is differentiable, except at $A_1$. This is easily remedied,
but, since it is irrelevant for the remainder, we omit the details.

%\vskip -0.2in
To continue, we first mention a curious result. The radius of curvature of the arcs $C_n$
is denoted by $r_n$ (see Figure \ref{fig:convexset1}). The radius of curvature by definition
equals $1/|z''(s)|$.
The following lemma can be shown by any persistent reader, since it uses only elementary planar
geometry. We omit the proof.

\begin{lemma} In the modified construction, $\lim_{n\rightarrow\infty}\,r_n$ exists
and is equal to $r_\infty=\dfrac{2\lambda}{1+\lambda}$.
\label{lem:radiusofcurv}
\end{lemma}

\begin{figure}[pbth]
\begin{center}
\includegraphics[height=4cm,width=5cm]{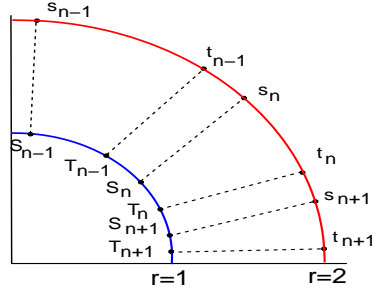}
\caption{The unit speed orbit $c(t)$ winds around the radius 2 circle. At times $t_n$, its projection
is $T_n$, and at times $s_n$, its projection is $S_n$.}
\label{fig:convexset2}
\end{center}
\end{figure}

Let $c(t)=2e^{it/2}$, the \emph{unit speed} curve along the radius 2 circle centered at the origin.
Near $t=0$, we see that $|\dot c_{\|}|=1$ and $|\dot c_{\perp}|=0$.
On the curve $c$, mark the times where $\Pi(t)=\Pi(c(t))$ equals $T_n$ by $t_n$, and times where
$\Pi(t)$ equals $S_n$ by $s_n$ (see Figure \ref{fig:convexset2}).
Again, using elementary planar geometry, the avid reader can prove the following lemma.

\begin{lemma} In the modified construction, the following limits exist:
\begin{equation*}
p_\infty :=\lim_{n\rightarrow \infty}\,\frac{t_{n-1}-s_{n}}{t_{n-1}} \quad \textrm{and} \quad
q_\infty :=\lim_{n\rightarrow \infty}\,\frac{s_{n}-t_{n}}{s_{n}} ,
\end{equation*}
and both are in $(0,1)$.
\label{lem:pn-and-qn}
\end{lemma}

\begin{theorem} The modified construction just described yields a differentiable curve $z(s)$
for the boundary of $K$ and $z'(s)$ is uniformly Lipschitz. The projection $\Pi$ onto this
set does \emph{not} have a directional derivative at $z=1$.
\label{thm:no-directional-deriv}
\end{theorem}

The remainder of this section is a sketch of the proof of Theorem \ref{thm:no-directional-deriv}.

\begin{proof}
We first establish that the curve is differentiable and its derivative is uniformly Lipschitz.
By construction, $\partial K$ consists of circular arcs $C_n$ (and a part of the unit circle),
and straight segments $T_{n-1}S_n$. So it is clear that $\partial K$ is $C^1$ and $C^2$ almost
everywhere. Where a segment and an arc
are glued together, the derivative of $z$ is a continuous function which is constant
on one side and has constant slope on the other. By Lemma \ref{lem:radiusofcurv},
this slope is bounded, and so the derivative is uniformly Lipschitz.

The only problematic point is $z=1$ where the arcs $C_n$ accumulate.
On the one hand, the line $\ell=\{1+it\,|\,t\in \mathbb{R}\}$ intersects $\partial K$ only at $z=1$.
The points $1$ and $T_n$ are contained in $K$. By convexity, the chords from $1$ to $T_n$ are
contained in $K$. The slope of these chords accumulate to the slope of $\ell$. Hence, $\ell$
is the tangent line at $z=1$. This establishes differentiability. With a little extra effort,
the argument in the previous paragraph also applies to $z=1$, establishing the Lipschitz condition there.

We next establish that the directional derivative of equation (\ref{eq:one-sidedetc}) does not exist near
$z=1$. If $t\in (t_{n},s_{n+1})$, then $\Pi(t)$ is part of a segment, and if $t\in(s_n,t_n)$, then by $\Pi(t)$ is in $C_n$.
Lemma \ref{lem:curvature} and Lemma \ref{lem:radiusofcurv} imply that for $n$ large (or
$t$ near 0)
\begin{equation}
\begin{array}{ccl}
\textrm{if } t\in (t_{n-1},s_{n}) &\textrm{ then }& |\Pi'(t)|=1 \quad \textrm{ and}\\
\textrm{if } t\in(s_n,t_n) &\textrm{ then }& |\Pi'(t)|=\dfrac{r_n}{1+r_n}
\rightarrow \dfrac{2\lambda}{1+3\lambda}<\dfrac 23.
\end{array}
\label{eq:piprime}
\end{equation}
This leads to the following equality:
\begin{equation}
\dfrac{\Pi(t_{n-1})-\Pi(0)}{t_{n-1}} =
\dfrac{\Pi(t_{n-1})-\Pi(s_{n})}{t_{n-1}-s_{n}} \,\dfrac{t_{n-1}-s_{n}}{t_{n-1}} +
\dfrac{\Pi(s_{n})-\Pi(0)}{s_{n}} \, \dfrac{s_{n}}{t_{n-1}} .
\label{eq:different speeds}
\end{equation}
Now if we assume that the one-sided directional derivative exists and insert equation (\ref{eq:piprime})
and $p_\infty$ defined in Lemma \ref{lem:pn-and-qn} into the previous equation, then, as $n$ tends
to infinity, we obtain
\begin{equation*}
\Pi_+'(0)= 1\cdot p_\infty + \Pi_+'(0)(1-p_\infty)  \quad \textrm{or} \quad
(\Pi_+'(0)-1) p_\infty=0.
\end{equation*}
It follows that $\Pi_+'(0)=1$ or that $p_\infty=0$. The latter is impossible by Lemma
\ref{lem:pn-and-qn}. Similarly, the equality
\begin{equation}
\dfrac{\Pi(s_{n})-\Pi(0)}{s_{n}} =
\dfrac{\Pi(s_{n})-\Pi(t_{n})}{s_{n}-t_{n}} \,\dfrac{s_{n}-t_{n}}{s_{n}} +
\dfrac{\Pi(t_{n})-\Pi(0)}{t_{n}} \, \dfrac{t_{n}}{s_{n}} ,
\label{eq:different speeds2}
\end{equation}
under the hypothesis that $\Pi_+'(0)$ exists, leads to the conclusion that $\Pi_+'(0)<\frac 23$.
Thus we have a contradiction and, hence, $\Pi_+'(0)$ cannot exist.
\end{proof}

There are two open questions related to this construction that are perhaps worth mentioning.
The first is: Can this construction be modified to obtain a curve $z(s)$ so that the
quantity $\Pi_+'$ fails to exist in a set of points that is dense in some open subset of $\Omega$?
The second question is: Is there a $C^2$ curve $\boldsymbol \sigma(u)$
in $\mathbb{R}^2$ for which the solution of equation (\ref{eq:ode3}) is not unique?
For such a curve, the one-sided derivative of the projection \emph{does} exist by
Lemma \ref{prop:Cp}. But the solution of the differential equation would not be unique.

%\vskip 0.4in\noindent
\section{Twice Differentiable Sets in $\mathbb{R}^2$.}
\label{chap:R2}

\noindent
For the next two sections, we return to the smooth (at least $C^2$) case.
It is interesting to translate the considerations in Section \ref{chap:C}
to the vector space $\mathbb{R}^n$. For simplicity, we start with $\mathbb{R}^2$.

Recall that $z:\mathbb{R}\rightarrow \mathbb{R}^2$ given by $\boldsymbol \sigma(u)$ is the unit speed
anti-clockwise parametrization of the boundary of the convex body.
Denote its unit tangent vector by ${\boldsymbol \sigma}'$, and $\bf n$ is the unit vector
normal to it, pointed into the body (to get a right-handed coordinate system
$({\boldsymbol \sigma}',\bf n)$). Note that $i z'$ in the complex notation becomes ${\bf n}$, and
that the complex $z''$ becomes $\kappa {\bf n}$. The sign of $\kappa$ --- nonnegative here ---
is determined by the convexity of the body. We now give
the translations of equations (\ref{eq:coords}) through (\ref{eq:ode}) to $\mathbb{R}^2$.

The position $c$ of a point outside the body is given by
\begin{equation}
c=\boldsymbol \sigma(u)-r{\bf n} .
\label{eq:coords2}
\end{equation}
(Note that $r$ is positive in $\Omega$.) Now let $c(t)$ be a smooth curve.
Differentiating with respect to $t$ while noting that
$\dot{\boldsymbol{n}}=-\kappa {\boldsymbol \sigma}'$, gives
\begin{equation}
\dot c = {\boldsymbol \sigma}'(1+r\kappa)\dot u - {\bf n}\dot r  .
\label{eq:cdot2}
\end{equation}
To write this in matrix form, define the matrix $\Sigma$ as having its first column equal
to ${\boldsymbol \sigma}'$ and its second equal to ${\bf n}$. Then
\begin{equation}
\dot c = \Sigma \begin{pmatrix}1+r\kappa & 0 \\ 0 & -1\end{pmatrix}
\begin{pmatrix}\dot u\\ \dot r\end{pmatrix} .
\label{eq:cdot3}
\end{equation}
Now $\Sigma$ is a unitary matrix (in fact, orthogonal) and so
\begin{equation}
\begin{pmatrix}\dot u\\ \dot r\end{pmatrix}= \begin{pmatrix}
(1+r\kappa)^{-1} & 0\\0 & -1\end{pmatrix}\, \Sigma^T \,\dot c .
\label{eq:ode2}
\end{equation}
Inverting the matrices gives us the next lemma.
\begin{lemma} Suppose $\boldsymbol{\sigma}$ is twice differentiable
and $c$ is a differentiable trajectory outside $K$. We have
\begin{equation}
\left\{\begin{array}{ccc}
\dot u &=& \dfrac{{\boldsymbol \sigma}'(u)\cdot \dot c}{1+r\kappa(u)}\\[0.35cm]
\dot r &=& -{\boldsymbol n}(u)\cdot \dot c
\end{array} \right. \quad ,
\label{eq:ode3}
\end{equation}
where $x\cdot y$ refers to the standard inner product in $\mathbb{R}^2$.
\label{lem:2dim}
\end{lemma}

In practice, the unit speed parametrization is not so useful, since for most curves
(including ellipses) explicit forms of such parametrizations are difficult or impossible
find \cite{ghrist}.
The method we use in Section \ref{chap:R3} will avoid this problem altogether.
For now, as a very simple illustration, consider the case where $K$ is the unit disk. Its circumference
is parametrized by ${\boldsymbol \sigma}(u)=(\cos u, \sin u)$. The radius of curvature
equals 1. Suppose $c=(x(t),y(t)$; then in this case the equations become
\begin{equation}
\left\{\begin{array}{ccc}
\dot u &=& \dfrac{-\dot x \sin u + \dot y \cos u}{1+r}\\[0.35cm]
\dot r &=& \dot x \cos u + \dot y \sin u
\end{array} \right. .
\end{equation}
We will look at more complicated examples in Section \ref{chap:examples}.

The assumption that $\boldsymbol \sigma(u)$ is $C^2$ is essential here. The example of
Section \ref{chap:counter-example} shows that $\dot u$ may not exist if the boundary is
$C^{1,1}$. The same holds in higher dimension (see Section \ref{chap:R3}),
though we will not repeat this observation.

%\vskip 0.4in\noindent
\section{Twice Differentiable Sets in $\mathbb{R}^3$.}
\label{chap:R3}

\noindent
In three dimensions, these equations become much more interesting. The
generalization from there to $\mathbb{R}^{n}$ is easy to guess;
one needs to replace the matrix $W$ below with the $n-1$ by $n-1$ Ricci
curvature tensor of the hypersurface. We will stick to $\mathbb{R}^3$ for simplicity.

Let ${\boldsymbol \sigma}(u,v)$ be a smooth (at least $C^2$) parametrization of a surface in $\mathbb{R}^3$
bounding a convex body $K$. Denote ${\partial_u\boldsymbol \sigma}$ by
$\boldsymbol \sigma_u$ and ${\partial_v\boldsymbol \sigma}$ by $\boldsymbol \sigma_v$.
To get a right-handed coordinate system
$({\boldsymbol \sigma}_u, {\boldsymbol \sigma}_v, {\bf n})$, we define the \emph{unit} normal ${\bf n}$:
\begin{displaymath}
{\bf n}=\dfrac{{\boldsymbol \sigma}_u\times {\boldsymbol \sigma}_v}
{{|\boldsymbol \sigma}_u\times {\boldsymbol \sigma}_v|} .
\end{displaymath}
To connect with the earlier sections, we will set things up in such a way that ${\bf n}$ points
\emph{into} the convex body $K$.
Let $\Sigma$ be the invertible $3\times 3$ matrix whose first column is $\boldsymbol \sigma_u$,
whose second column is $\boldsymbol \sigma_v$, and whose third is ${\bf n}$.
Note that the above assumptions imply that the determinant of this matrix is positive.
Furthermore, let $W$ be the $2 \times 2$ \emph{Weingarten matrix}, that is (see \cite{DoCarmo}),
\begin{equation}
\left\{\begin{array}{ccc}
{\bf n}_u &=& a {\boldsymbol \sigma}_u+ b{\boldsymbol \sigma}_v\\[0.35cm]
{\bf n}_v &=& c {\boldsymbol \sigma}_u+ d{\boldsymbol \sigma}_v\\[0.35cm]
\end{array} \right.
\quad \textrm{where} \quad
\begin{pmatrix} a&c\\b&d\end{pmatrix}=W .
\label{eq:weingartenmatrix}
\end{equation}
The eigenvalues of $-W$ are the principal curvatures of the surface (see \cite{DoCarmo}).
By convexity --- recall that ${\bf n}$ points away from $\Omega$ --- these eigenvalues are $\emph{nonnegative}$.

\begin{theorem} Suppose the parametrization $\boldsymbol \sigma$ of the boundary of $K$ is twice differentiable
and $c$ is a differentiable trajectory outside $K$. Let $I$ be the $2\times 2$ identity matrix. We have
\begin{equation}
\begin{pmatrix}\dot u\\ \dot v \\ \dot r\end{pmatrix}=
\begin{pmatrix} (I - rW)^{-1} & \vline &0\\
\hline 0 & \vline & -1 \end{pmatrix} \, \Sigma^{-1} \, \dot c .
\label{eq:diffeq}
\end{equation}
\label{thm:3dim}
\end{theorem}

\vskip -0.2in\noindent
\begin{proof} Coordinatize $\Omega$, the space outside the convex body, as follows (see equation
(\ref{eq:coords2})):
\begin{displaymath}
P:=(u,v,r)\rightarrow {\boldsymbol \sigma}(u,v)-r{\bf n} ,
\end{displaymath}
where $r>0$ in $\Omega$. Let $c(t):=(u(t),v(t),r(t))$ be a smooth curve in $\Omega$, then
\begin{displaymath}
\begin{array}{ccl}
\dot c &=& {\boldsymbol \sigma}_u \dot u + {\boldsymbol \sigma}_v \dot v -
{\bf n} \dot r - r {\bf n}_u \dot u - r {\bf n}_v \dot v \\[0.2cm]
&= & {\boldsymbol \sigma}_u((1-ra)\dot u-rc\dot v)+
{\boldsymbol \sigma}_v(-rb \dot u+(1-rd)\dot v) - {\bf n} \dot r ,
\end{array}
\end{displaymath}
where we have used equation (\ref{eq:weingartenmatrix}). In matrix form, this gives
\begin{displaymath}
\dot c = \Sigma \, \begin{pmatrix} I - rW & \vline &0\\
\hline 0 & \vline & -1 \end{pmatrix}
\begin{pmatrix}\dot u\\ \dot v \\ \dot r\end{pmatrix}.
\end{displaymath}
Since the eigenvalues of $-W$ are nonnegative and $r$
is positive in $\Omega$, this relation can be inverted to yield the theorem.
\end{proof}

Given the vector $\dot c$ of the trajectory of the observer, the derivative
of the projection is given by ${\boldsymbol \sigma}_u \dot u+{\boldsymbol \sigma}_v \dot v$.
Thus equation (\ref{eq:diffeq}) implies the following.

\begin{corollary} If the surface in $\mathbb{R}^3$ is twice differentiable, then
the directional derivative exists.
\label{cor:deriv-proj}
\end{corollary}

For existence and uniqueness of the solution of this differential equation, the conditions
of Theorem \ref{thm:3dim} are insufficient. We need to require that
the right-hand side is Lipschitz in $(u,v,r)$ and continuous in $t$ (see \cite{HSD}).

\begin{corollary} (Local) existence and uniqueness of the solution of differential equation
(\ref{eq:diffeq}) is guaranteed if $W$ is Lipschitz and $c$ is continuously differentiable.
\label{cor:exist-uniq}
\end{corollary}

\vskip .0in \noindent
\begin{proof} The continuity in $t$ of the right-hand side of equation (\ref{eq:diffeq}) is guaranteed
by the smoothness of $c$. We now check the Lipschitz condition in $(u,v,r)$. Lipschitz in $r$ is straightforward.
Since ${\boldsymbol \sigma}$ is $C^2$,  $\Sigma$ is differentiable and invertible, and therefore
so is $\Sigma^{-1}$. Hence, $\Sigma^{-1}$ is certainly Lipschitz. The question then is whether $(I-rW)^{-1}$
is Lipschitz. The following computation shows that if $W$ is Lipschitz, then so is $(I-rW)^{-1}$.

From the previous discussion, we are allowed to hold $r$ constant.
Let us denote the operator norm by $\|\cdot\|$ and set $x=(u,v)$. Then
\begin{equation*}
\begin{array}{l}
 \|(I-rW(x_2))^{-1}-(I-rW(x_1))^{-1}\| \\
\hskip 0.5in =  |(I-rW(x_2))^{-1}\left[r(W(x_2)-W(x_1)\right](I-rW(x_1))^{-1}| \\
\hskip 0.5in \leq   r\|(I-rW(x_2))^{-1}\| \| W(x_2)-W(x_1)\| \|(I-rW(x_1))^{-1}\|  .
\end{array}
\end{equation*}
The first and last term are at most 1 by convexity (the eigenvalues of $W$ are nonpositive).
The middle term is Lipschitz by assumption.
\end{proof}

If we want the solution of equation (\ref{eq:diffeq}) to depend differentiably on time and
initial conditions, we need to require even more.

\begin{proposition} Let $\phi(t,u_0,v_0,r_0)$ be the solution of equation (\ref{eq:diffeq}) with
initial condition $(u_0,v_0,r_0)$ at $t=0$.  Let $p\geq 1$ an integer. If $c(t)$ is $C^{p+1}$ and
${\boldsymbol \sigma}(u,v)$ is $C^{p+2}$, then $\phi$ is $C^p$ in $(t,u_0,v_0,r_0)$.
\label{prop:Cp}
\end{proposition}

\vskip .0in \noindent
\begin{proof} In \cite[Chapter 32]{arnold}, it is proved that if $\dot x=f(t,x)$ and $f$ is
$C^p$ for some integer $p\geq 1$, then the initial value problem is also $C^p$.
Inspection of equation (\ref{eq:diffeq}) shows that the conditions of the proposition are
sufficient to guarantee that its right-hand side is $C^p$.
\end{proof}

%\vskip 0.4in\noindent
\section{Differentiability of the Distance Function.}
\label{chap:dist-diffb}

\noindent
If we restrict ourselves to the distance function, we can do much better for \emph{its} regularity than the
results in Section \ref{chap:R3}. In this section, we assume that $K$ is closed and convex, but make
no further regularity assumption.

\begin{lemma} Suppose that $x_i$ is in $\Omega$ and $\lim x_i=x_\infty$. Let $y_i\in\Pi(x_i)$. Then
any convergent subsequence of $\{y_i\}$ converges to a point $y_\infty$ in $\Pi(x_\infty)$.
\label{lem:conts}
\end{lemma}

\vskip -0.0in\noindent
\begin{proof} Suppose the opposite; then pick a subsequence of the $i$ so that $\{y_i\}$ is
convergent to a point $y_\infty$ that is not in $\Pi(x_\infty)$. Let $y\in\Pi(x_\infty)$; then
for some positive $\Delta$, we have $d(x_\infty,y_\infty)=d(x_\infty,y)+\Delta$.
Now take $\epsilon$ arbitrarily small. Using the triangle inequality, we see that for $n$ large enough
$d(x_\infty,x_n)+d(x_n,y_n)+d(y_n,y_\infty)\geq d(x_\infty,y_\infty)$, and so
\begin{displaymath}
d(x_n,y_n)\geq d(x_\infty,y_\infty)-d(x_\infty,x_n)-d(y_n,y_\infty)\geq d(x_\infty,y)+\Delta-2\epsilon .
\end{displaymath}
On the other hand, using the triangle inequality for large $n$ again,
\begin{displaymath}
d(x_n,y)\leq d(x_\infty,y)+d(x_n,x_\infty)<d(x_\infty,y)+\epsilon .
\end{displaymath}
Taking the two together, we get that $d(x_n,y_n)> d(x_n,y)$, which contradicts the fact that $y_n$ minimizes the distance from $x_n$ to $K$.
\end{proof}

If $K$ is convex, then $\Pi(x_\infty)$ is a single point $y$. Thus every convergent subsequence
of $\{y_i\}$ in the Lemma \ref{lem:conts} must converge to that point. Since the $\{y_i\}$ are confined
to a bounded portion of a closed set, they always have a convergent subsequence.
Thus we have proved a nice corollary:

\begin{corollary} The projection $\Pi:\Omega\rightarrow K$ is continuous.
\label{cor:conts}
\end{corollary}

For the remainder of this section, denote by $\boldsymbol{\theta}$ the unit vector in the
direction of $x-\Pi(x)$.

\begin{lemma} For any $x\in \Omega$ and any vector ${\bf v}$, we have
\begin{displaymath}
d(x+t{\bf v},\Pi(x))-d(x,\Pi(x))= t\boldsymbol{\theta} \cdot{\bf v}+o(t) .
\end{displaymath}
\label{lem:plaut}
\end{lemma}

\vskip -0.4in\noindent
\begin{proof} The difference between the squares of the distances can be written using the
standard inner product of $\mathbb{R}^n$:
\begin{displaymath}
d(x+t {\bf v},\Pi(x))^2-d(x,\Pi(x))^2 = 2t[x-\Pi(x)]\cdot{\bf v} + t^2 |{\bf v}|^2 .
\end{displaymath}
The left-hand side factors, so that we can divide by $d(x+t{\bf v},\Pi(x))+d(x,\Pi(x))\sim 2|x-\Pi(x)|$,
from which the lemma follows easily.
\end{proof}

\vskip -0.2in\noindent
\begin{lemma} For any $x\in \Omega$ and any vector ${\bf v}$, we have
\begin{displaymath}
\begin{array}{l}
d(x+t{\bf v},\Pi(x+t{\bf v}))-d(x,\Pi(x+t{\bf v})) \\
\hskip 0.5in =t\boldsymbol{\theta} \cdot{\bf v}+ t\dfrac{\left[\Pi(x)-\Pi(x+t{\bf v})\right]}{|x-\Pi(x)|}\cdot{\bf v} +o(t) .
\end{array}
\end{displaymath}
\label{lem:plautplus}
\end{lemma}

\vskip -0.2in\noindent
\begin{proof} The computation is essentially the same as that in the proof of Lemma \ref{lem:plaut}.
\begin{equation}
\begin{array}{l}
\hskip -0.3in
d(x+t {\bf v},\Pi(x+t{\bf v}))^2-d(x,\Pi(x+t{\bf v}))^2 \\[0.1cm]
\hskip 0.3in = 2t[x-\Pi(x+t{\bf v})]\cdot{\bf v} + t^2 |{\bf v}|^2 \\[0.1cm]
\hskip 0.3in = 2t[x-\Pi(x)]\cdot{\bf v} + 2t[\Pi(x)-\Pi(x+t{\bf v})]\cdot{\bf v}+t^2 |{\bf v}|^2 .
\end{array}
\end{equation}
As before, we divide by $d(x+t{\bf v},\Pi(x))+d(x,\Pi(x))\sim 2|x-\Pi(x)|$.
\end{proof}

In our case where $K$ is convex, Corollary \ref{cor:conts} implies the second term in the right-hand
side of the statement of Lemma \ref{lem:plautplus} is also $o(t)$. This is important in the following
theorem. (An altogether different proof of Theorem \ref{thm:differentiable} can be found in
\cite[Chapter 2]{giaquinta}.)

\vskip 0.0in\noindent
\begin{theorem} For any $x\in \Omega$ and any vector ${\bf v}$, we have
\begin{displaymath}
d(x+t{\bf v},\Pi(x+t{\bf v}))-d(x,\Pi(x))= t\boldsymbol{\theta} \cdot{\bf v}+o(t) .
\end{displaymath}
In other words, the distance function to a convex set in $\mathbb{R}^n$ is differentiable on $\Omega$.
\label{thm:differentiable}
\end{theorem}

\vskip 0.0in\noindent
\begin{proof} We note that by definition of the projection $\Pi$, we have
\begin{displaymath}
d(x+t {\bf v},\Pi(x+t {\bf v}))-d(x,\Pi(x))\leq d(x+t {\bf v},\Pi(x))-d(x,\Pi(x)) .
\end{displaymath}
And thus by Lemma \ref{lem:plaut},
\begin{displaymath}
d(x+t {\bf v},\Pi(x+t {\bf v}))-d(x,\Pi(x))\leq t{\boldsymbol \theta}\cdot{\bf v}+o(t) .
\end{displaymath}
Similarly, if we use Lemma \ref{lem:plautplus} and Corollary \ref{cor:conts} to make sure
that $t|\Pi(x)-\Pi(x+t{\bf v})|=o(t)$, we obtain the reverse inequality
\begin{displaymath}
d(x,\Pi(x))\leq d(x,\Pi(x+t {\bf v}))=d(x+t {\bf v},\Pi(x+t {\bf v}))-
t{\boldsymbol \theta} \cdot{\bf v}+o(t) .
\end{displaymath}
\end{proof}

It turns out Lemma \ref{lem:plaut} actually holds not just in $\mathbb{R}^n$ but in any
Riemannian manifold (see \cite{Plaut}). The proof we gave of Lemma \ref{lem:plautplus} is valid
only in $\mathbb{R}^n$. The proof of Theorem \ref{thm:differentiable} is a modified version of
the proof that appeared in \cite{foote}. However, Theorem \ref{thm:differentiable} also holds for Alexandrov spaces
(a generalization of Riemannian manifolds) with nonpositive or nonnegative curvature, according
to \cite[Exercise 4.5.11]{BBI}.

%\vskip 0.4in\noindent
\section{Examples.}
\label{chap:examples}

\begin{figure}[pbth]
\begin{center}
\includegraphics[width=4.5cm]{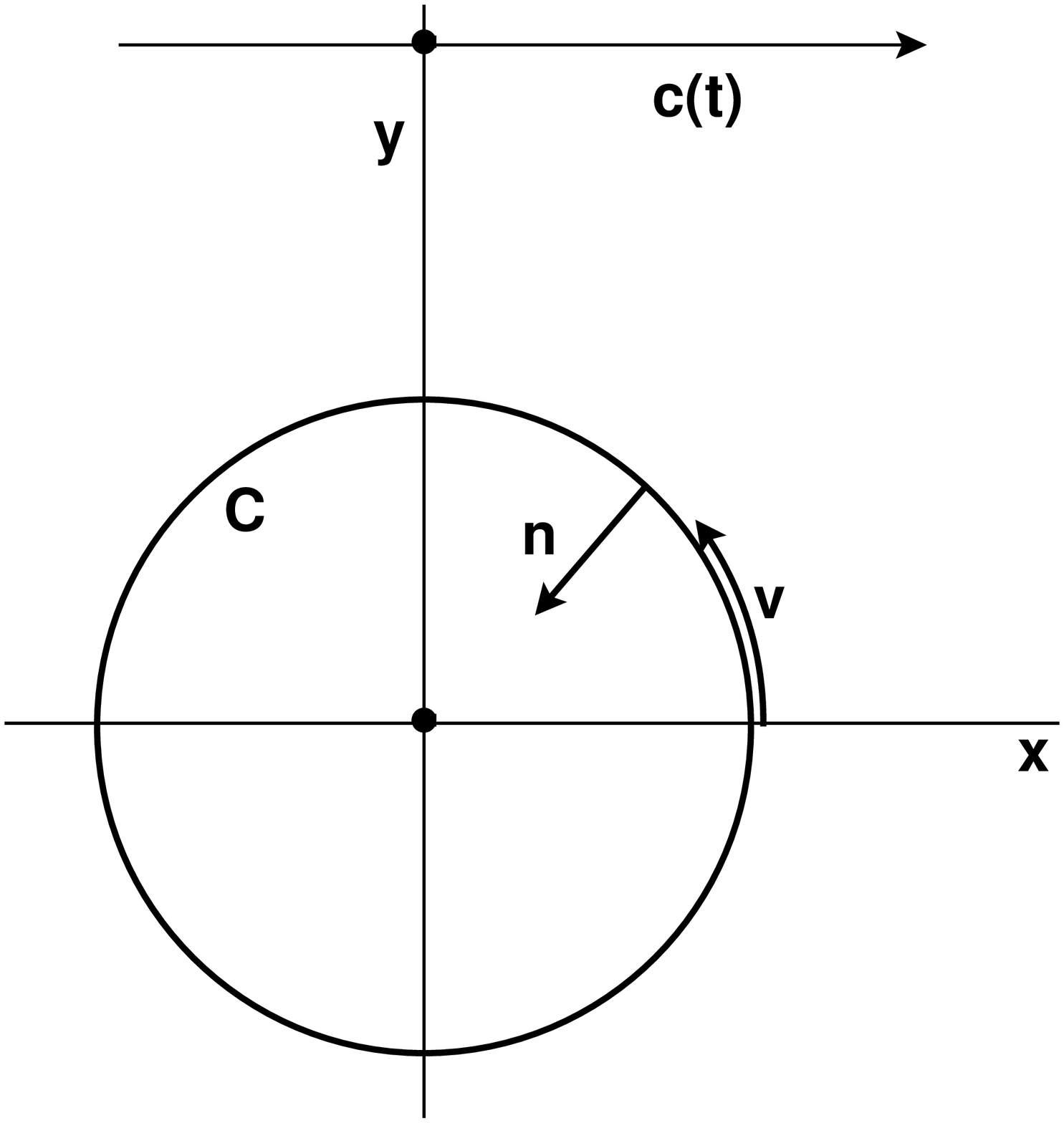} \quad \quad \quad
\includegraphics[width=4.5cm]{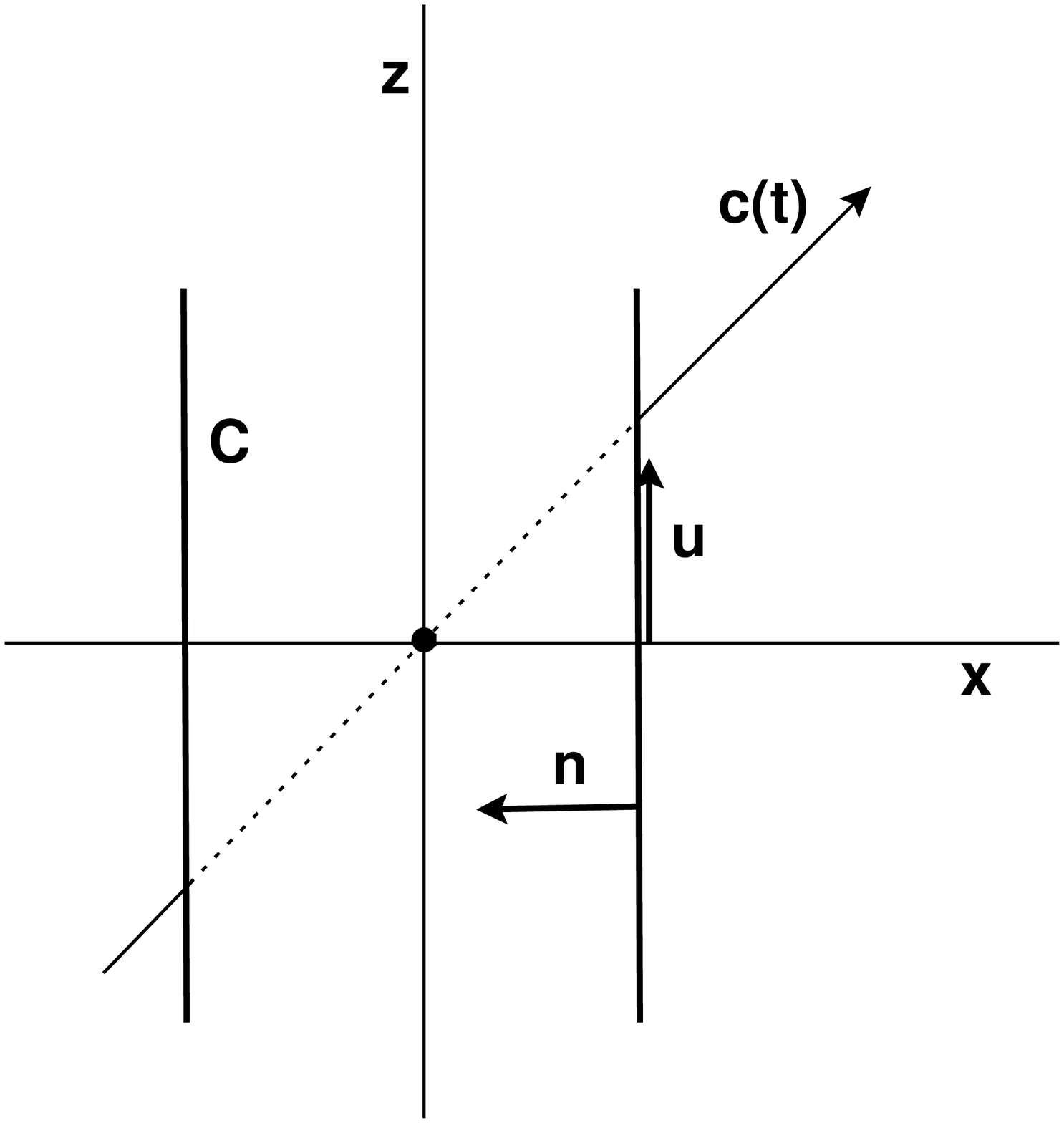}
\caption{On the left, the projection onto the $xy$-plane of the cylinder $C$ and the path $c$.
On the right, the projection onto the $xz$-plane.}
\label{fig:cylinder}
\end{center}
\end{figure}

\noindent
We close with a simple illustration of the equations discussed in Section \ref{chap:R3}. The
situation is sketched in Figure \ref{fig:cylinder}. This example also illustrates the fact
that straight lines in $\Omega$ do not generally project to geodesics on the surface.

Let $\sigma$ be the embedding of the cylinder of radius $\kappa^{-1}$ in $\mathbb{R}^3$. We choose
a parametrization so that its orientation is consistent with that of the previous sections.
\begin{equation*}
{\boldsymbol \sigma}:(u,v)\rightarrow (\kappa^{-1}\cos v, \kappa^{-1}\sin v, u) .
\end{equation*}
For positive constants $a$ and $b$ let the trajectory $c(t)$ be given by
\begin{equation*}
c(t)=(t, b+\kappa^{-1}, at) .
\end{equation*}

%\vskip -0.2in
In the notation of Section \ref{chap:R3}, we calculate the matrix $\Sigma$ (whose first column is
$\boldsymbol \sigma_u$, the second, $\boldsymbol \sigma_v$, and the third, ${\bf n}$).
\begin{equation*}
\Sigma= \begin{pmatrix} 0&-\kappa^{-1}\sin v&-\cos v \\
                        0&\kappa^{-1}\cos v &-\sin v\\
                        1 & 0          & 0\end{pmatrix} \quad \Longrightarrow \quad
\Sigma^{-1}= \begin{pmatrix} 0         & 0        & 1\\
                             -\kappa \sin v & \kappa \cos v &0\\
                             -\cos v     &  -\sin v      &0 \end{pmatrix} .
\end{equation*}
The Weingarten matrix $W$ is computed as in equation (\ref{eq:weingartenmatrix}):
\begin{equation*}
W= \begin{pmatrix} 0&0\\
                   0&-\kappa\end{pmatrix} \quad \Longrightarrow \quad
(I-rW)^{-1}= \begin{pmatrix} 1&0\\
                   0&(1+\kappa r)^{-1}\end{pmatrix} .
\end{equation*}
The projection of $c(t)$ onto the cylinder together with its distance to the cylinder is thus
given by the solution of equation (\ref{eq:diffeq}), together with initial conditions. Noting that
$\dot c=(1,0,a)$, these equations simplify to
\begin{equation*}
\begin{array}{ccl}
\dot u & = & a  \\[0.1cm]
\dot v & = & \dfrac{-\kappa\sin v}{1+\kappa r} \\[0.1cm]
\dot r & = & \cos v
\end{array} \quad \textrm{where} \quad
\begin{array}{ccl}
u(0) & = & 0  \\[0.1cm]
v(0) & = & \frac{\pi}{2} \\[0.1cm]
r(0) & = & b
\end{array} .
\end{equation*}

\begin{figure}[pbth]
\begin{center}
\includegraphics[width=4.5cm]{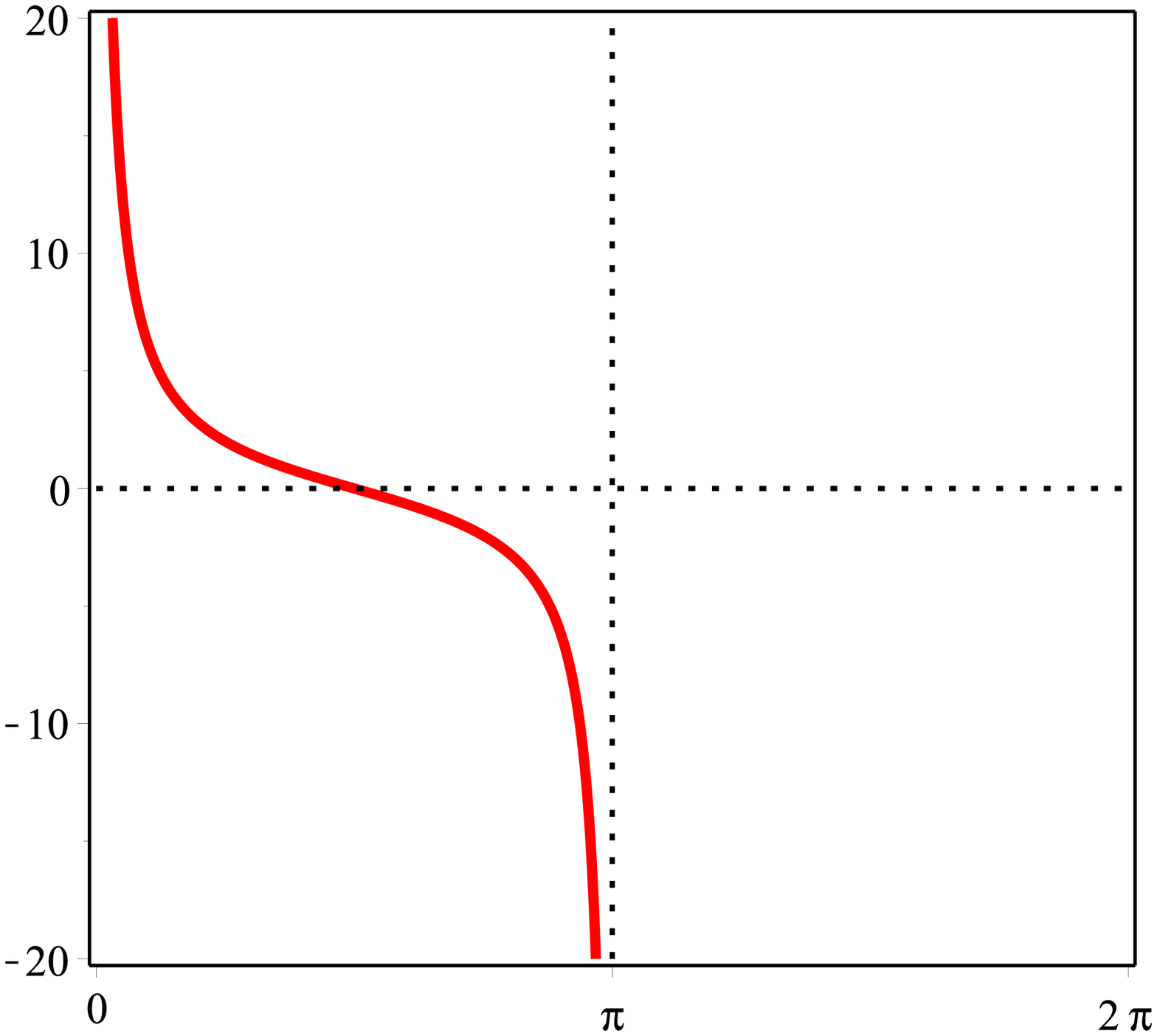} \quad \quad \quad
\includegraphics[width=4.5cm]{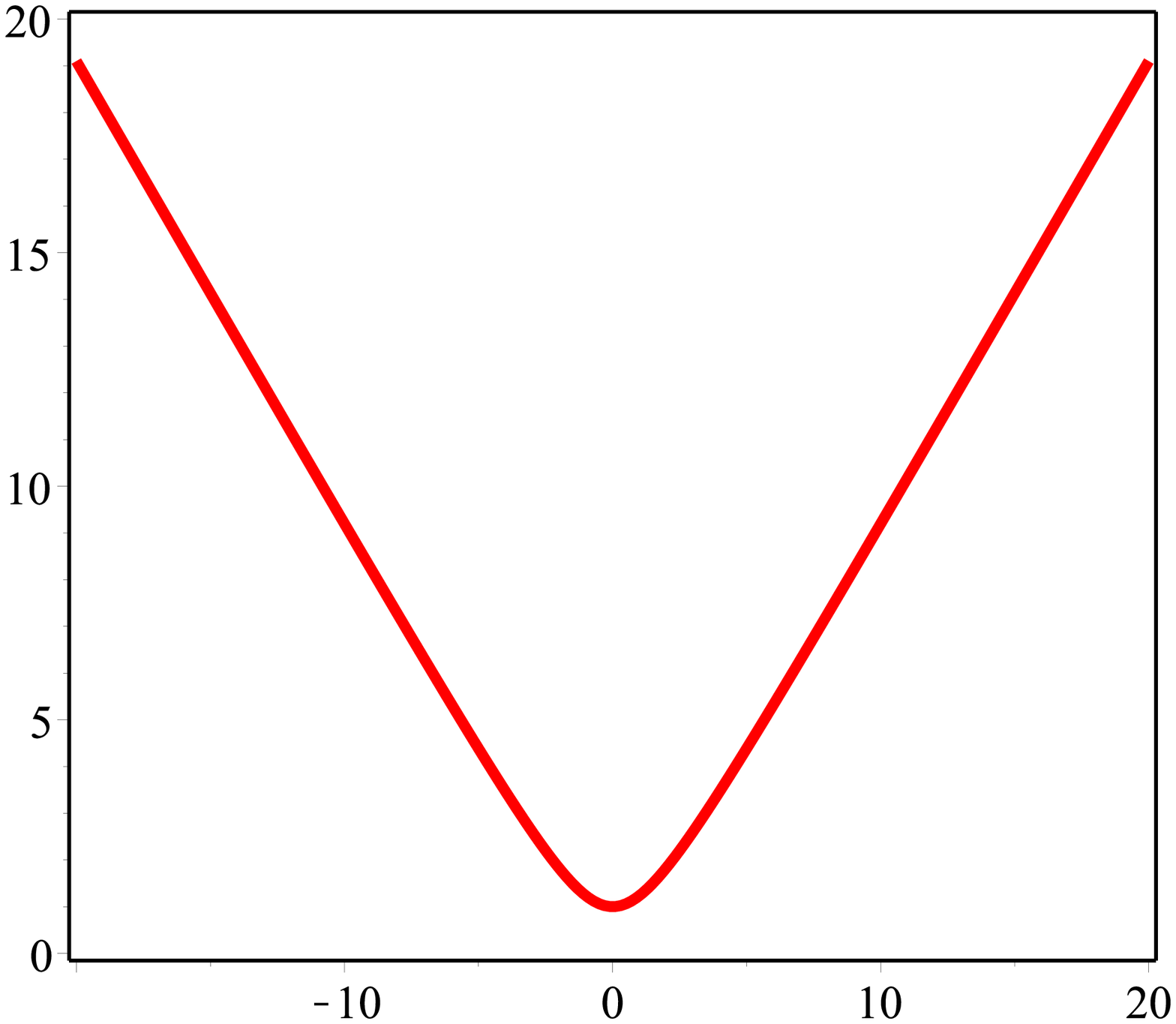}
\caption{On the left, the projection of the trajectory $C:t\rightarrow (t,2,t)$ onto the
cylinder $C:(u,v)\rightarrow (\cos v, \sin v, u)$. Horizontal is $v$,
vertical is $u$. On the right, the distance $r(t)$ between $c(t)$ and its projection onto the cylinder $C$.
Horizontal is $t$, vertical is $r$.}
\label{fig:cylinder3-4}
\end{center}
\end{figure}

This somewhat obscure-looking nonlinear system nonetheless has a simple solution.
Due to the translational symmetry, the segment from $c(t)$ to its projection lies in the
plane parallel to the $xy$-plane. Thus the distance $r(t)$ and the angle $v(t)$ can be found
by inspection of Figure \ref{fig:cylinder} on the left. We get
\begin{equation*}
\begin{array}{ccl}
u(t) & = & at  \\[0.1cm]
v(t) & = & \arccos \dfrac{t}{\sqrt{\ell^2+t^2}} \\[0.1cm]
r(t) & = & \sqrt{\ell^2+t^2}-\kappa^{-1}
\end{array}  ,
\end{equation*}
where we have set $\ell \equiv \kappa^{-1}+b$. By inverting the second of these, we obtain its inverse $t(v)=\ell (\tan v)^{-1}=- \ell \tan(v+\frac{\pi}{2})$. In the $vu$-plane, this solution is therefore
a reparametrization of the curve $(v,-a\ell \tan(v+\frac{\pi}{2}))$. We have drawn this curve
and the curve $(t,r(t))$ for $t \in [-20,20]$ in Figure \ref{fig:cylinder3-4}.

\begin{figure}[pbth]
\begin{center}
\includegraphics[width=6.5cm,height=4.0cm]{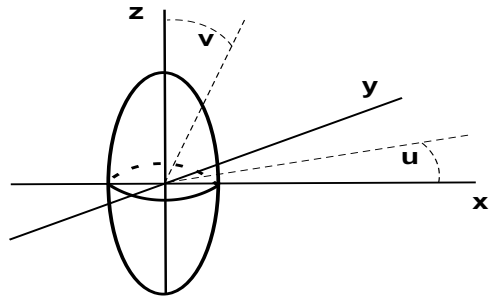}
\caption{An ellipsoid of revolution in spherical coordinates. Here, $u$ is the angle with the
axis of the ``north pole" and $v$ is the angle in the $xy$ plane with the positive $x$ axis.}
\label{fig:ellipsoid1}
\end{center}
\end{figure}

We now turn to a slightly more challenging example.
This time $\boldsymbol \sigma$ is the embedding of the ellipsoid of revolution in $\mathbb{R}^3$.
We choose a parametrization so that its orientation is consistent with that of the previous sections.
This includes making sure that the normal is pointing into the surface.
\begin{equation*}
{\boldsymbol \sigma}:(u,v)\rightarrow (\cos u \sin v, \sin u \sin v, \kappa^{-1/2} \cos v) .
\end{equation*}
Thus when $\kappa\in(0,1)$, the ellipsoid looks a little like a cigar, as in Figure \ref{fig:ellipsoid1},
and when  $\kappa\in(1,\infty)$, the ellipsoid is a flattened one.

With the same conventions as before, and following exactly the same procedure,
we compute the matrix $S$ of tangent vectors and the Weingarten matrix $W$:
\begin{equation*}
\Sigma= \begin{pmatrix} -\sin u \sin v & \cos u \cos v & {\dfrac {-\cos u \sin v }{\sqrt {1+(\kappa-1)
\cos^{2} v}}}\\[0.4cm]
\cos u \sin v & \sin u \cos v & {\dfrac {-\sin u \sin v }{\sqrt {1+(\kappa-1)  \cos^{2} v}}}\\[0.4cm]
0 & \dfrac{-\sin v}{\sqrt{\kappa}} & {\dfrac {-\sqrt{\kappa}\cos v }{\sqrt {1+(\kappa-1)  \cos^{2} v}}}
\end{pmatrix}  ;
\end{equation*}
\begin{equation*}
W= \begin{pmatrix} \dfrac{-1}{\left(1+(\kappa-1)  \cos^{2} v\right)^{1/2}} & 0\\
                   0 & \dfrac{-\kappa}{\left(1+(\kappa-1)  \cos^{2} v\right)^{3/2}}
                   \end{pmatrix} .
\end{equation*}
After setting
\begin{equation*}
\begin{array}{ccl}
q_1 & = & \sqrt{1+(\kappa-1)\cos^2 v}
\end{array} ,
\end{equation*}
and some algebra, the equations of motion (\ref{eq:diffeq}) become:
\begin{equation}
\begin{array}{ccl}
\dot u & = & \dfrac{q_1}{(r+q_1)\sin v} \;(-\dot x\,\sin u+\dot y \, \cos u) \\[0.4cm]
\dot v & = & \dfrac{ q_1\, \kappa}{r\kappa+q_1^3}
\;(\dot x \, \cos u \cos v+\dot y \,\sin u \cos v-\dot z\,\kappa^{-1/2} \sin v)\\[0.4cm]
\dot r & = & \dfrac{1}{q_1} \;(\dot x \cos u \sin v+ \dot y \sin u \sin v+ \dot z \sqrt{\kappa} \cos v)
\end{array} .
\label{eq:ellipsoidexample}
\end{equation}
Here $c(t)=(x(t),y(t),z(t))$ is the a priori given trajectory.

\begin{figure}[pbth]
\begin{center}
\includegraphics[width=7.3cm]{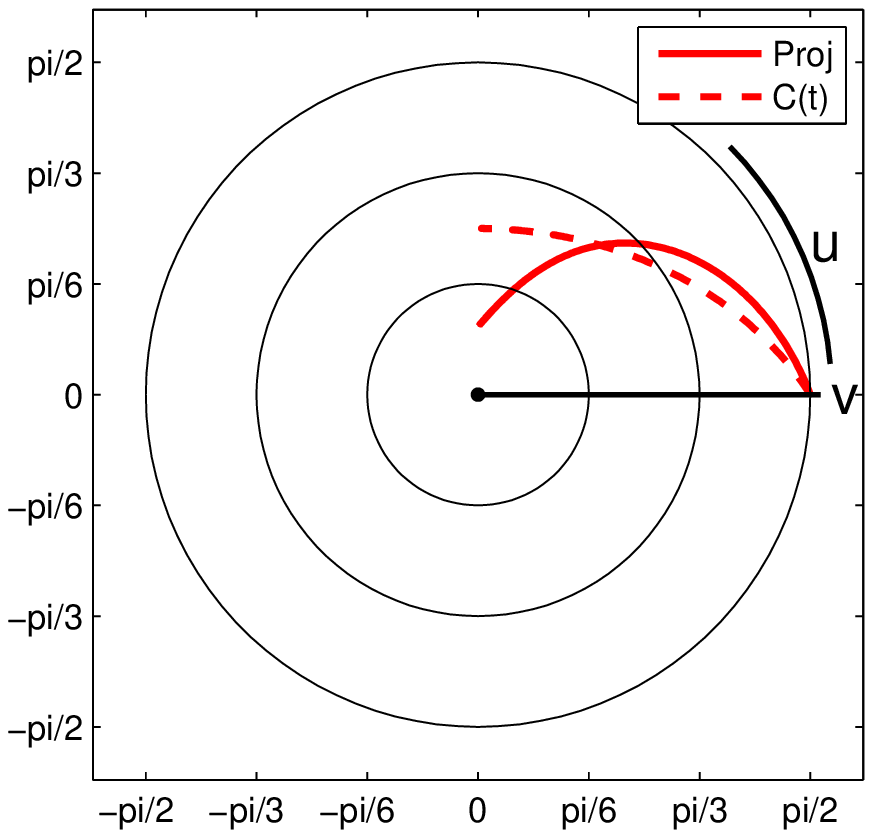}
\includegraphics[width=5.2cm]{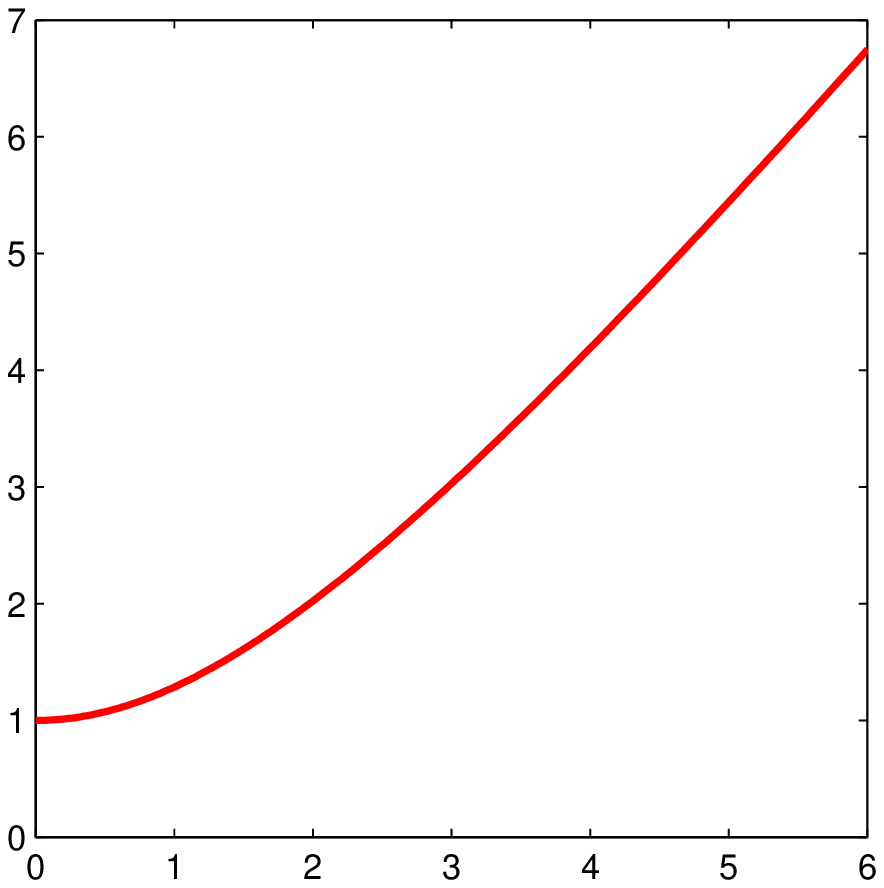}
\caption{In the figure on the left, the projection $\Pi(c(t))$ (continuous curve)
and the trajectory $c(t)$ (dashed curve) are displayed for $t\in[0,100]$. We use $(u,v)$ coordinates
(see text). In the figure on the right, $r(t)$ is displayed for $t\in[0,6]$ ($r$ is vertical, time
is horizontal).}
\label{fig:ellipsoid2}
\end{center}
\end{figure}

As an example, we set $\kappa=1/9$ and $c(t)=(2,t,t)$ and numerically solved equation
(\ref{eq:ellipsoidexample}) setting the initial condition $(u(0),v(0),r(0))=(0,\pi/2,1)$ and
$(\dot x, \dot y, \dot z)=(0,1,1)$. The results are displayed in Figure
\ref{fig:ellipsoid2}. The figure on the left shows the orbit of the projection $\Pi(c(t))$
for $t\in[0,100]$ (the continuous curve) and the orbit $c(t)$ (the dashed curve).
We display only the $(u,v)$ coordinates; the radial coordinate is $v$, so that the
circles show where $v$ is constant (namely 0, $\frac{\pi}{6}$, $\frac{2\,\pi}{6}$, and
$\frac{3\,\pi}{6}$) and the angular coordinate is $u$ ($u$ equals the angle with the positive $x$-axis
indicated in the figure).
The figure on the right shows the function $r(t)$. For clarity, we only show the function for $t\in[0,6]$.
These figures were produced with the Matlab routine ode45 with AbsTol and RelTol set to $10^{-9}$.

To be absolutely certain that the result is reliable, we recall that we should be able to recover
$c(t)=(2,t,t)$ from the solution of equation (\ref{eq:ellipsoidexample}) as follows (see Section
\ref{chap:R3}).
\begin{displaymath}
c(t)={\boldsymbol \sigma}(u(t),v(t))-r(t){\bf n}(t) .
\end{displaymath}
For the displayed curve we found
\begin{displaymath}
\max_{t\in[0,100]}\,|{\boldsymbol \sigma}(u(t),v(t))-r(t){\bf n}(t)-c(t)|\leq 3\cdot 10^{-7}.
\end{displaymath}

\begin{acknowledgment}{Acknowledgment.}
I am indebted to John Milnor for a very useful conversation which led to equation \ref{eq:cdot},
and to Robert Lyons for the numerical results presented in Figure \ref{fig:ellipsoid2}.
\end{acknowledgment}

%\vskip 0.4 in

\begin{biog}
\item[J. J. P. Veerman] received his Ph.D. from Cornell University. He has held
visiting positions in the U.S. (Rockefeller University, Stony Brook University, Georgia Tech, Penn State),
as well as in Spain, Brazil, Italy, and Greece. He is currently at Portland State
University in Oregon, USA, where he is Professor of Mathematics and Affiliate Professor
of Physics.
\begin{affil}
Maseeh Department of Mathematics and Statistics, Portland State University, Portland, OR 97201.\\
veerman@pdx.edu
\end{affil}
\end{biog}

\vspace{\fill}
\end{document}